\documentclass[11pt]{amsart}
\usepackage{amsfonts}
\usepackage{amsmath}
\usepackage{amssymb}

\theoremstyle{plain}

\newtheorem{remark}{Remark}

\numberwithin{equation}{section}

\input{tcilatex}

\begin{document}
\author[S.~H. Saker]{S. ~H. Saker}
\address{Department of Mathematics, Faculty of Science, Mansoura University\\
Mansoura, 35516, Egypt}
\email{shsaker@mans.edu.eg}
\title[Hardy Type Inequalities]{New Hardy-Type Inequalities Via Opial
Inequalities}
\maketitle

\begin{abstract}
In this paper, we will prove several new inequalities of Hardy's type with
explicit constants. The main results will be proved by making use of some
generalizations of Opial's type inequalities and H\"{o}lder's inequality. To
the best of the author's knowledge this method has not been used before in
investigation of this type of inequalities. From these inequalities one can
establish some new inequalities of differential forms which are inequalities
of Wirtinger's type.
\end{abstract}

\section{Introduction}

The classical Hardy inequality (see \cite{H}) states that for $f\geq 0$ and
integrable over any finite interval $(0,x)$ and $f^{p}$ is integrable and
convergent over $(0,\infty )$ and $p>1,$ then 
\begin{equation}
\int_{0}^{\infty }\left( \frac{1}{x}\left( \int_{0}^{x}f(t)dt\right)
dx\right) ^{p}\leq \left( \frac{p}{p-1}\right) ^{p}\int_{0}^{\infty
}f^{p}(x)dx,  \label{h0}
\end{equation}%
unless $f=0.$ The constant $\left( p/\left( p-1\right) \right) ^{p}$ is the
best possible. This inequality has been proved by Hardy in 1925, but it has
been appeared as the continuous version of a discrete inequality in his work
in 1920 when he aimed to find a new elementary proof of Hilbert's inequality
for double series. The discrete version of (\ref{h0}) is given by (see \cite%
{Notes}) the discrete inequality%
\begin{equation*}
\sum_{n=1}^{\infty }\left( \frac{1}{n}\sum_{k=1}^{n}a_{k}\right) ^{p}\leq
\left( \frac{p}{p-1}\right) \sum_{n=1}^{\infty }a_{n}\text{, \ (}a_{n}>0%
\text{, \ \ }p>1).
\end{equation*}%
During the past decades, the study of Hardy inequalities (continuous and
discrete) or Hardy operators focused on the investigations of new
inequalities or operators with weighted functions. These results are of
interest and important in analysis not only because the mappings are optimal
in the sense that the size of weight classes cannot be improved, but also
because the weight conditions themselves are interest. This intensively
investigated area of mathematical analysis resulted in the publication of
numerous research papers and books, we refer the reader to the books \cite%
{K1,K2,K3} and the papers \cite{Og, BH, Conv, Ka}. The inequality (\ref{h0})
was extended to the form 
\begin{equation}
\int_{a}^{b}\left( \left( \int_{a}^{x}f(t)dt\right) ^{q}u(x)dx\right) ^{%
\frac{1}{q}}\leq C\left( \int_{a}^{b}f^{p}(x)v(x)dx\right) ^{\frac{1}{q}},
\label{h1}
\end{equation}%
with $a,$ $b$ real numbers, satisfying $a<b$, and $u,$ $v$ are positive and
measurable functions in the interval $(a,b)$, $p,$ $q$ are positive
parameters satisfying $0<q<\infty $ and $1\leq p<\infty .$ The main idea in
this type of inequalities is to find the optimal value of the constant $C$
(see \cite{K1}). Another type of inequalities of Hardy's type has been
proved by Beesack in \cite{BH}. In particular, Beesack proved some
inequalities of the form 
\begin{equation}
\int_{a}^{b}s(x)\left( \int_{0}^{x}f(t)dt\right) ^{p}dx\leq
\int_{a}^{b}r(x)f^{p}(x)dx,\text{ }  \label{BH}
\end{equation}%
where the weighted functions $r$ and $s$ satisfy the Euler-Lagrange
differential equation%
\begin{equation}
\frac{d}{dx}(r(x)(y^{^{\prime }})^{p-1})+s(x)y^{p-1}(x)=0.  \label{EL}
\end{equation}%
The method of the proofs in \cite{BH} depends on the existence of positive
solution $y$ of (\ref{EL}) which satisfies $y>0$ and $y^{^{\prime }}>0$ on
the interval $(a,b).$ The idea of the present paper has been appeared during
my work in some papers deal with disconjugacy and gaps between zeros of some
differential equations (\cite{S1, S2, S4}), which require some inequalities
of the form 
\begin{equation}
\int_{a}^{b}s(x)y^{p}dx\leq C\int_{a}^{b}r(x)(y^{^{\prime }}(x))^{p}dx,
\label{W}
\end{equation}%
where $C$ is the constant of the inequality depends on the weighted
functions and needs to be determined. The differential forms of Hardy's
inequalities (like \ref{W}) are of Wirtinger's type (see \cite{ABOS, S3}).
This type of inequalities are useful in the analysis of qualitative behavior
of solutions of differential equations and can be used in investigations of
disconjugacy and disfocality of solutions (see \cite{S1, S2, S4}).

\bigskip

Kaijser et al. \cite{Ka} pointed out that the inequality (\ref{h0}) is just
a special case of the much more general (Hardy-Knopp type) inequality 
\begin{equation}
\int_{0}^{\infty }\Phi \left( \frac{1}{x}\int_{0}^{x}f(t)dt\right) \frac{dx}{%
x}\leq \int_{0}^{\infty }\Phi \left( f(x)\right) \frac{dx}{x},  \label{Conv}
\end{equation}%
where $\Phi $ is a convex function on $(0,\infty ).$ Recently the inequality
(\ref{Conv}) has been generalized by Kaijser et al. \cite{Conv} by using the
properties of convex functions and by employed Jensen's inequality and
Fubini's Theorem.

\bigskip

Our aim in this paper is to prove some inequalities with weighted functions
of Hardy's type by making use of some generalizations of Opial inequalities
and H\"{o}lder's inequality. The method that, we will apply in this paper is
different from the techniques used in the above mentioned papers and books.
Throughout the paper some special cases\ are derived. The differential forms
of these inequalities which are inequalities of Wirtinger's type can be
obtained and due to the limited space the details will be left to the
interested reader.

\section{Main Results}

In this section, we will prove the main results by making use of the H\"{o}%
lder inequality\textit{\ }%
\begin{equation}
\int_{a}^{b}\left\vert f(x)g(x)\right\vert dx\leq \left[ \int_{a}^{b}\left%
\vert f(x)\right\vert ^{p}dx\right] ^{\frac{1}{p}}\left[ \int_{a}^{b}\left%
\vert g(x)\right\vert ^{q}dx\right] ^{\frac{1}{q}},  \label{H0}
\end{equation}%
where $a$, $b\in \mathbb{I}$ and $f;\,g$\ $\in C(\mathbb{I},\mathbb{R}),$ $%
p>1$\ such that $\frac{1}{p}+\frac{1}{q}=1,$ and some generalizations of
Opial's inequality will be stated latter. The Opial inequality states that
(see \cite{opial}): If $y$ is absolutely continuous on $[a,b]$ and $%
y(a)=y(b)=0,$ then 
\begin{equation}
\int_{a}^{b}\left\vert y(x)\right\vert \left\vert y^{^{\prime
}}(x)\right\vert dx\leq \frac{\left( b-a\right) }{4}\int_{a}^{b}\left\vert
y^{^{\prime }}(x)\right\vert ^{2}dx,  \label{opial}
\end{equation}%
with a best constant $1/4.$ Since the discovery of Opial inequality much
work has been done, and many papers which deal with new proofs, various
generalizations and extensions have appeared in the literature. In further
simplifying the proof of the Opial inequality which had already been
simplified by Olech \cite{Olech}, Beesack \cite{B1}, Levinson \cite{L1},
Mallows \cite{M1} and Pederson \cite{P1}, it is proved that if $y$ is real
absolutely continuous on $(a,b)$ and with $y(a)=0,$ then 
\begin{equation}
\int_{a}^{b}\left\vert y(x)\right\vert \left\vert y^{^{\prime
}}(x)\right\vert dx\leq \frac{b}{2}\int_{a}^{b}\left\vert y^{^{\prime
}}(x)\right\vert ^{2}dx.  \label{B1}
\end{equation}%
These inequalities and their extensions and generalizations are the most
important and fundamental inequalities in the analysis of qualitative
properties of solutions of different types of differential equations, (see 
\cite{S1, S2, S4}).

Throughout the paper, all functions are assumed to be positive and
measurable and all the integrals appear in the inequalities are exist and
finite.

\smallskip

\textbf{Theorem 2.1.} \textit{Assume that }$r$\textit{, }$s$\textit{\ be
positive and continuous functions on the interval }$(a,b).$\textit{\ Then }%
\begin{equation}
\left( \int_{a}^{b}r(x)\left( \int_{a}^{x}f(t)dt\right) dx\right) ^{2}\leq
\int_{a}^{b}\frac{R^{2}(x,b)}{s(x)}dx\left( \int_{a}^{b}s(x)\left(
f(x)\right) ^{2}dx\right) ,  \label{H11}
\end{equation}%
\textit{for all functions}$f\geq 0$\textit{, where}%
\begin{equation*}
R(x,b)=\int_{x}^{b}r(x)dx.
\end{equation*}%
\textbf{Proof.} Let $F(x)=\int_{a}^{x}f(t)dt.$ It is clear that $F(a)=0$ and 
$F^{^{\prime }}(x)=f(x).$ This gives us that 
\begin{equation*}
\int_{a}^{b}r(x)\left( \int_{a}^{x}f(t)dt\right) dx=\int_{a}^{b}r(x)F(x)dx.
\end{equation*}%
Integrating by parts the right hand side, we have 
\begin{equation*}
\int_{a}^{b}r(x)\left( \int_{a}^{x}f(t)dt\right) dx=-\left.
R(x,b)F(x)\right\vert _{a}^{b}+\int_{a}^{b}R(x,b)F^{^{\prime }}(x)dx.
\end{equation*}%
Using the assumptions $R(b,b)=0$ and $F(a)=0,$ we have 
\begin{eqnarray*}
\int_{a}^{b}r(x)\left( \int_{a}^{x}f(t)dt\right) dx
&=&\int_{a}^{b}R(x,b)F^{^{\prime }}(x)dx \\
&=&\int_{a}^{b}\frac{R(x,b)}{\sqrt{s(x)}}\sqrt{s(x)}F^{^{\prime }}(x)dx.
\end{eqnarray*}%
Applying the inequality (\ref{H0}) with $p=q=2$, we see that 
\begin{eqnarray*}
\int_{a}^{b}r(x)\left( \int_{a}^{x}f(t)dt\right) dx &\leq &\left(
\int_{a}^{b}\frac{R^{2}(x,b)}{s(x)}dx\right) ^{\frac{1}{2}}\left(
\int_{a}^{b}s(x)\left( F^{^{\prime }}(x)\right) ^{2}dx\right) ^{\frac{1}{2}}
\\
&=&\left( \int_{a}^{b}\frac{R^{2}(x,b)}{s(x)}dx\right) ^{\frac{1}{2}}\left(
\int_{a}^{b}s(x)\left( f(x)\right) ^{2}dx\right) ^{\frac{1}{2}},
\end{eqnarray*}%
which is the desired inequality (\ref{H11}). The proof is complete.

\smallskip

As in the proof of Theorem 2.1, by putting $F(x)=\int_{x}^{b}f(t)dt,$ we can
can prove the following result.

\smallskip

\textbf{Theorem 2.2.} \textit{Assume that }$r$\textit{, }$s$\textit{\ be
positive and continuous functions on the interval }$(a,b).$\textit{\ Then }%
\begin{equation}
\left( \int_{a}^{b}r(x)\left( \int_{x}^{b}f(t)dt\right) dx\right) ^{2}\leq
\int_{a}^{b}\frac{R^{2}(a,x)}{s(x)}dx\left( \int_{a}^{b}s(x)\left(
f(x)\right) ^{2}dx\right) ,  \label{H12}
\end{equation}%
\textit{for all functions}$f\geq 0$\textit{, where}%
\begin{equation*}
R(a,x)=\int_{a}^{x}r(x)dx.
\end{equation*}

In the following, we apply the inequality (\ref{B1}) to prove a new
inequality of Hardy's type.

\smallskip

\textbf{Theorem 2.3}\textit{. Assume that }$r$\textit{\ be a positive and
continuous function on the interval }$(a,b).$\textit{\ Then }%
\begin{equation}
\int_{a}^{b}r(x)\left( \int_{a}^{x}f(t)dt\right) ^{2}dx\leq
(b-a)\sup_{a<x<b}\left( \int_{x}^{b}r(x)dx\right) \int_{a}^{b}\left(
f(x)\right) ^{2}dx,  \label{Hl0}
\end{equation}%
\textit{for all functions }$f\geq 0$\textit{.}

\textbf{Proof.} We proceed as in the proof of Theorem 2.1 to get 
\begin{equation*}
\int_{a}^{b}r(x)\left( \int_{a}^{x}f(t)dt\right)
^{2}dx=\int_{a}^{b}r(x)F^{2}(x)dx.
\end{equation*}%
Integrating by parts the right hand side, we have 
\begin{equation*}
\int_{a}^{b}r(x)\left( \int_{a}^{x}f(t)dt\right) ^{2}dx=-\left.
R(x,b)F^{2}(x)\right\vert _{a}^{b}+2\int_{a}^{b}R(x,b)F(x)F^{^{\prime
}}(x)dx.
\end{equation*}%
Using the assumptions $R(b,b)=0$ and $F(a)=0,$ we have 
\begin{eqnarray*}
\int_{a}^{b}r(x)\left( \int_{a}^{x}f(t)dt\right) ^{2}dx
&=&2\int_{a}^{b}R(x,b)F(x)F^{^{\prime }}(x)dx \\
&\leq &2\sup_{a<x<b}R(x,b)\int_{a}^{b}F(x)F^{^{\prime }}(x)dx.
\end{eqnarray*}%
Applying the inequality (\ref{B1}), since $F(a)=0$, we obtain 
\begin{eqnarray*}
\int_{a}^{b}r(x)\left( \int_{a}^{x}f(t)dt\right) dx &\leq
&(b-a)\sup_{a<x<b}R(x,b)\int_{a}^{b}\left( F^{^{\prime }}(x)\right) ^{2}dx \\
&=&(b-a)\sup_{a<x<b}R(x,b)\left( \int_{a}^{b}\left( f(x)\right)
^{2}dx\right) ,
\end{eqnarray*}%
which is the desired inequality (\ref{Hl0}). The proof is complete.

\smallskip

As in the proof of Theorem 2.3, one can prove the following result.

\smallskip

\textbf{Theorem 2.4.}\textit{\ Assume that }$r$\textit{\ be positive and
continuous function on the interval }$(a,b).$\textit{\ Then }%
\begin{equation*}
\int_{a}^{b}r(x)\left( \int_{x}^{b}f(t)dt\right) ^{2}dx\leq
(b-a)\sup_{a<x<b}\left( \int_{a}^{x}r(x)dx\right) \int_{a}^{b}\left(
f(x)\right) ^{2}dx,
\end{equation*}%
\textit{for all functions }$f\geq 0$.

\smallskip

In the following, we will apply a generalization of (\ref{B1}) due to
Beesack \cite{B1} to prove a new inequality of Hardy's type. The inequality
due to Beesack states that: If $y~$\ is an absolutely continuous function on 
$[a,b]$ with $y(a)=0$ (or $y(b)=0$) then 
\begin{equation}
\int_{a}^{b}\left\vert y(t)\right\vert \left\vert y^{^{\prime
}}(t)\right\vert dt\leq \frac{1}{2}\int_{a}^{b}\frac{1}{r(t)}%
dt\int_{a}^{b}r(t)\left\vert y^{^{\prime }}(t)\right\vert ^{2}dt,  \label{B2}
\end{equation}%
where $r(t)$ is positive and continuous function with $\int_{a}^{b}dt/r(t)<%
\infty .$ Using this inequality on the term%
\begin{equation*}
2\int_{a}^{b}R(x,b)F(x)F^{^{\prime }}(x)dx\leq
2\sup_{a<x<b}R(x,b)\int_{a}^{b}F(x)F^{^{\prime }}(x)dx,
\end{equation*}%
we see that 
\begin{eqnarray*}
2\int_{a}^{b}R(x,b)F(x)F^{^{\prime }}(x)dx &\leq &2\sup_{a<x<b}R(x,b)\frac{1%
}{2}\left( \int_{a}^{b}\frac{dx}{s(x)}\right) \int_{a}^{b}s(x)\left(
F^{^{\prime }}(x)\right) ^{2}dx \\
&=&\sup_{a<x<b}R(x,b)\left( \int_{a}^{b}\frac{dx}{s(x)}\right)
\int_{a}^{b}s(x)\left( f(x)\right) ^{2}dx.
\end{eqnarray*}%
Using this inequality and proceeding as in the proof of Theorem 2.3, we have
the following inequality.

\smallskip

\textbf{Theorem 2.5.}\textit{\ Assume that }$r$, $s$\textit{\ be positive
and continuous functions in the interval }$(a,b).$\textit{\ Then }%
\begin{equation*}
\int_{a}^{b}r(x)\left( \int_{a}^{x}f(t)dt\right) ^{2}dx\leq
\sup_{a<x<b}R(x,b)\left( \int_{a}^{b}\frac{dx}{s(x)}\right)
\int_{a}^{b}s(x)\left( f(x)\right) ^{2}dx,
\end{equation*}%
\textit{for all functions }$f\geq 0$.

\smallskip

\textbf{Theorem 2.6.}\textit{\ Assume that }$r$\textit{, }$s$\textit{\ be
positive and continuous functions on the interval }$(a,b).$\textit{\ Then }%
\begin{equation*}
\int_{a}^{b}r(x)\left( \int_{x}^{b}f(t)dt\right) ^{2}dx\leq
\sup_{a<x<b}R(a,x)\left( \int_{a}^{b}\frac{dx}{s(x)}\right)
\int_{a}^{b}s(x)\left( f(x)\right) ^{2}dx,
\end{equation*}%
\textit{for all functions }$f\geq 0$.

\smallskip

In the following, we apply the Maroni inequality \cite{M2} which states
that: If $y~$\ is an absolutely continuous function on $[a,b]$ with $y(a)=($%
or $y(b)=0),$ then%
\begin{equation}
\int_{a}^{b}\left\vert y(t)\right\vert \left\vert y^{^{\prime
}}(t)\right\vert dt\leq \frac{1}{2}\left( \int_{a}^{b}\left( \frac{1}{r(t)}%
\right) ^{p-1}dt\right) ^{\frac{2}{p}}\left( \int_{a}^{b}r(t)\left\vert
y^{^{\prime }}(t)\right\vert ^{q}dt\right) ^{\frac{2}{q}},  \label{M1}
\end{equation}%
where $\int_{a}^{b}\left( 1/r(t)\right) ^{p-1}dt<\infty $, $p\geq 1$ and $%
\frac{1}{p}+\frac{1}{q}=1.$

\smallskip

Applying the inequality (\ref{M1}) on the term $\sup_{a<x<b}R(x,b)%
\int_{a}^{b}F(x)F^{^{\prime }}(x)dx$, we see that 
\begin{eqnarray*}
&&2\int_{a}^{b}R(x,b)F(x)F^{^{\prime }}(x)dx\leq
2\sup_{a<x<b}R(x,b)\int_{a}^{b}F(x)F^{^{\prime }}(x)dx \\
&\leq &\sup_{a<x<b}R(x,b)\left( \int_{a}^{b}\left( \frac{1}{s(x)}\right)
^{p-1}dx\right) ^{\frac{2}{p}}\left( \int_{a}^{b}s(x)\left( F^{^{\prime
}}(x)\right) ^{q}dx\right) ^{\frac{2}{q}} \\
&=&\sup_{a<x<b}R(x,b)\left( \int_{a}^{b}\left( \frac{1}{s(x)}\right)
^{p-1}dx\right) ^{\frac{2}{p}}\left( \int_{a}^{b}s(x)\left( f(x)\right)
^{q}dx\right) ^{\frac{2}{q}}.
\end{eqnarray*}%
Using this inequality and proceed as in the proof of Theorem 2.3, we have
the following results.

\smallskip

\textbf{Theorem 2.7. }\textit{Assume that }$r,s$\textit{\ be positive and
continuous functions in the interval }$(a,b)$\textit{and }$p>1$\textit{\
such that }$\frac{1}{p}+\frac{1}{q}=1.$\textit{\ Then }%
\begin{eqnarray*}
\int_{a}^{b}r(x)\left( \int_{a}^{x}f(t)dt\right) ^{2}dx &\leq
&\sup_{a<x<b}R(x,b)\left( \int_{a}^{b}\left( \frac{1}{s(x)}\right)
^{p-1}dx\right) ^{\frac{2}{p}} \\
&&\times \left( \int_{a}^{b}s(x)\left( f(x)\right) ^{q}dx\right) ^{\frac{2}{q%
}},
\end{eqnarray*}%
\textit{for all functions }$f\geq 0$.

\smallskip

\textbf{Theorem 2.8. }\textit{Assume that }$r$\textit{\ be a positive
function in the interval }$(a,b)$\textit{\ such that }$\int_{a}^{b}\left(
1/R(a,x)\right) ^{p-1}dx<\infty $\textit{, }$p\geq 1$\textit{\ and }$\frac{1%
}{p}+\frac{1}{q}=1.$\textit{\ Then }%
\begin{eqnarray*}
\int_{a}^{b}r(x)\left( \int_{x}^{b}f(t)dt\right) ^{2}dx &\leq
&\sup_{a<x<b}R(a,x)\left( \int_{a}^{b}\left( \frac{1}{s(x)}\right)
^{p-1}dx\right) ^{\frac{2}{p}} \\
&&\times \left( \int_{a}^{b}s(x)\left( f(x)\right) ^{q}dx\right) ^{\frac{2}{q%
}}
\end{eqnarray*}%
\textit{for all functions }$f\geq 0.$

\smallskip

Yang \cite{Y1} simplified the Beesack proof and extended the inequality (\ref%
{B2}) and proved that:\ If $y$ is an absolutely continuous function on $[a,b]
$ with $y(a)=0$ (or $y(b)=0)$ then 
\begin{equation}
\int_{a}^{b}q(t)\left\vert y(t)\right\vert \left\vert y^{^{\prime
}}(t)\right\vert dt\leq \frac{1}{2}\int_{a}^{b}\frac{1}{r(t)}%
dt\int_{a}^{b}r(t)q(t)\left\vert y^{^{\prime }}(t)\right\vert ^{2}dt,
\label{Y}
\end{equation}%
where $r(t)$ is a positive and continuous function with $%
\int_{a}^{b}dt/r(t)<\infty $ and $q(t)$ is a positive, bounded and
non-increasing function on $[a,b].$ Applying the inequality (\ref{Y}) on the
term $\int_{a}^{b}R(x,b)F(x)F^{^{\prime }}(x)dx$, we see that 
\begin{eqnarray*}
&&\int_{a}^{b}R(x,b)F(x)F^{^{\prime }}(x)dx \\
&\leq &\frac{1}{2}\left( \int_{a}^{b}\frac{1}{s(x)}dx\right)
\int_{a}^{b}R(x,b)s(x)\left( F^{^{\prime }}(x)\right) ^{2}dx \\
&=&\frac{1}{2}\left( \int_{a}^{b}\frac{1}{s(x)}dx\right)
\int_{a}^{b}R(x,b)s(x)\left( f(x)\right) ^{2}dx.
\end{eqnarray*}%
This gives us the following results.

\textbf{Theorem 2.9. }\textit{Assume that }$r$\textit{, }$s$\textit{\ be
positive functions in the interval }$(a,b)$\textit{\ such that }$%
\int_{a}^{b}dt/s(t)<\infty $\textit{\ and }$R(x,b)$\textit{\ is a positive,
bounded and non-increasing function on }$[a,b].\,$\textit{Then }%
\begin{equation*}
\int_{a}^{b}r(x)\left( \int_{a}^{x}f(t)dt\right) ^{2}dx\leq \left(
\int_{a}^{b}\frac{dx}{s(x)}\right) \int_{a}^{b}R(x,b)s(x)\left( f(x)\right)
^{2}dx,
\end{equation*}%
\textit{for all functions }$f\geq 0.$\textit{\ }

\smallskip

\textbf{Theorem 2.10.}\textit{\ Assume that }$r$\textit{, }$s$\textit{\ be
positive and continuous functions on the interval }$(a,b)$\textit{\ such
that }$\int_{a}^{b}dt/s(t)<\infty $\textit{\ and }$R(a,x)$\textit{\ is a
positive, bounded and non-increasing function on }$[a,b].\,$\textit{Then }%
\begin{equation*}
\int_{a}^{b}r(x)\left( \int_{x}^{b}f(t)dt\right) ^{2}dx\leq \left(
\int_{a}^{b}\frac{dx}{s(x)}\right) \left( \int_{a}^{b}R(a,x)s(x)\left(
f(x)\right) ^{2}dx\right) ,
\end{equation*}%
\textit{for all functions }$f\geq 0.$

\smallskip

In the following, we will apply the inequality due to Hua \cite{H1} to prove
a new inequality of Hardy's type. The inequality due to Hua states that: If $%
y~$is an absolutely continuous function with $y(a)=0$ (or $y(b)=0)$, then 
\begin{equation}
\int_{a}^{b}\left\vert y(t)\right\vert ^{p}\left\vert y^{^{\prime
}}(t)\right\vert dt\leq \frac{(b-a)^{p}}{p+1}\int_{a}^{b}\left\vert
y^{^{\prime }}(t)\right\vert ^{p+1}dt,  \label{H1}
\end{equation}%
where $p$ is a positive integer.

\smallskip

\textbf{Theorem 2.11.}\textit{\ Assume that }$r$\textit{\ be positive and
continuous function on the interval }$(a,b)$\textit{\ and }$p$\textit{\ is a
positive integer. Then }%
\begin{equation}
\int_{a}^{b}r(x)\left( \int_{a}^{x}f(t)dt\right) ^{p+1}dx\leq
(b-a)^{p}\sup_{a<x<b}\left( \int_{x}^{b}r(t)dt\right) \int_{a}^{b}\left(
f(x)\right) ^{p+1}dx,  \label{Hl1}
\end{equation}%
\textit{\ for all functions }$f\geq 0.$

\textbf{Proof.} Let $F(x)=\int_{a}^{x}f(t)dt.$ It is clear that $F(a)=0$ and 
$F^{^{\prime }}(x)=f(x)>0.$ This gives that%
\begin{equation*}
\int_{a}^{b}r(x)\left( \int_{a}^{x}f(t)dt\right)
^{p+1}dx=\int_{a}^{b}r(x)F^{p+1}(x)dx.
\end{equation*}%
Integrating by parts the left hand side, we have 
\begin{eqnarray*}
\int_{a}^{b}r(x)\left( \int_{a}^{x}f(t)dt\right) ^{p+1}dx &=&-\left.
R(x,b)F^{p+1}(x)\right\vert _{a}^{b} \\
&&+\left( p+1\right) \int_{a}^{b}R(x,b)F^{p}(x)F^{^{\prime }}(x)dx.
\end{eqnarray*}%
Using the assumptions $R(b,b)=0$ and $F(a)=0,$ we have 
\begin{eqnarray}
\int_{a}^{b}r(x)\left( \int_{a}^{x}f(t)dt\right) ^{p+1}dx
&=&(p+1)\int_{a}^{b}R(x,b)F^{p}(x)F^{^{\prime }}(x)dx  \notag \\
&\leq &(p+1)\sup_{a<x<b}R(x,b)\int_{a}^{b}F^{p}(x)F^{^{\prime }}(x)dx.
\label{c1}
\end{eqnarray}%
Applying the inequality (\ref{H1}), since $F(a)=0$, on the term $%
\int_{a}^{b}F^{p}(x)F^{^{\prime }}(x)dx$, we have%
\begin{equation}
\int_{a}^{b}F^{p}(x)F^{^{\prime }}(x)dx\leq \frac{(b-a)^{p}}{p+1}%
\int_{a}^{b}\left( F^{^{\prime }}(x)\right) ^{p+1}dx.  \label{c2}
\end{equation}%
Substituting (\ref{c1}) into (\ref{c2}), we have that 
\begin{equation*}
\int_{a}^{b}r(x)\left( \int_{a}^{x}f(t)dt\right) ^{p+1}dx\leq
(b-a)^{p}\sup_{a<x<b}R(x,b)\int_{a}^{b}\left( f(x)\right) ^{p+1}dx,
\end{equation*}%
which is the desired inequality (\ref{Hl1}). The proof is complete.

\smallskip

As in the proof of Theorem 2.11 one can prove the following result.

\smallskip

\textbf{Theorem 2.12.} \textit{Assume that }$r$\textit{\ be positive and
continuous function on the interval }$(a,b)$\textit{\ and }$p$\textit{\ is a
positive integer. Then }%
\begin{equation*}
\int_{a}^{b}r(x)\left( \int_{x}^{b}f(t)dt\right) ^{p+1}dx\leq
(b-a)^{p}\sup_{a<x<b}\left( \int_{a}^{x}r(t)dt\right) \int_{a}^{b}\left(
f(x)\right) ^{p+1}dx,
\end{equation*}%
\textit{for all functions }$f\geq 0.$

\smallskip

Boyd and Wong \cite{BW} extended the inequality (\ref{H1}) for general
values of $p>0$ and proved that if $y~$is an absolutely continuous\ function
on $[a,b]$ with $y(a)=0$ (or $y(b)=0)$, then 
\begin{equation}
\int_{a}^{b}s(t)\left\vert y(t)\right\vert ^{p}\left\vert y^{^{\prime
}}(t)\right\vert dt\leq \frac{1}{\lambda _{0}(p+1)}\int_{a}^{b}r(t)\left%
\vert y^{^{\prime }}(t)\right\vert ^{p+1}dt,  \label{BW1}
\end{equation}%
where $r$ and $s$ are nonnegative functions in $C^{1}[a,b]$, $\lambda _{0}$
is the smallest eigenvalue of the boundary value problem%
\begin{equation*}
(r(t)\left( u^{^{\prime }}(t)\right) ^{p})^{^{\prime }}=\lambda s^{^{\prime
}}(t)u^{p}(t),
\end{equation*}%
with $u(a)=0$ and $u(b)=0$ for which $u^{^{\prime }}>0$ in $[a,b].$ Applying
the inequality (\ref{BW1}) on the term $(p+1)\int_{a}^{b}R(x,b)F^{p}(x)F^{^{%
\prime }}(x)dx,$ we have 
\begin{equation}
(p+1)\int_{a}^{b}R(x,b)F^{p}(x)F^{^{\prime }}(x)dx\leq \frac{1}{\lambda _{0}}%
\int_{a}^{b}s(t)(F^{^{\prime }}(x))^{p+1}dx,  \label{R1}
\end{equation}%
where $r$ and $s$ are nonnegative functions in $C^{1}[a,b]$, and $\lambda
_{0}$ is the smallest eigenvalue of the boundary value problem%
\begin{equation}
(R(x,b)\left( u^{^{\prime }}(x)\right) ^{p})^{^{\prime }}=\lambda
s^{^{\prime }}(x)u^{p}(x),  \label{Rb}
\end{equation}%
where $R(x,b)=\int_{x}^{b}r(t)dt$, with $u(a)=0$ and $u(b)=0$ for which $%
u^{^{\prime }}>0$ in $[a,b].$ Using (\ref{R1}) and proceeding as in the
proof of Theorem 2.11, we have the following result.

\smallskip

\textbf{Theorem 2.13.}\textit{\ Assume that }$r,$\textit{\ }$s$\textit{\ be
positive and continuous functions on }$(a,b)$\textit{\ and }$p$\textit{\ is
a positive integer. Then }%
\begin{equation*}
\int_{a}^{b}r(x)\left( \int_{a}^{x}f(t)dt\right) ^{p+1}dx\leq \frac{1}{%
\lambda _{0}}\int_{a}^{b}s(t)\left( f(x)\right) ^{p+1}dx,
\end{equation*}%
\textit{for all functions }$f\geq 0$\textit{\ where }$\lambda _{0}$\textit{\
is the smallest eigenvalue of the boundary value problem (\ref{Rb})}$.$

\smallskip

In the following, we will apply the Calvert inequality (see \cite{Cal}) 
\begin{equation}
\int_{a}^{b}\left\vert y(t)\right\vert ^{^{p}}\left\vert y^{^{\prime
}}(t)\right\vert dt\leq \frac{1}{p+1}\left( \int_{a}^{b}s^{\frac{-1}{p}%
}(t)bt\right) ^{p}\int_{a}^{b}s(t)\left\vert y^{^{\prime }}(t)\right\vert
^{p+1}dt,  \label{ag}
\end{equation}%
where $y(a)=0$ (or $y(b)=0)$ and $s(t)>0$, to prove a new inequality of
Hardy's type. Applying the inequality (\ref{ag}) on the term $%
\int_{a}^{b}F^{p}(x)F^{^{\prime }}(x)dx$, we have 
\begin{equation*}
\int_{a}^{b}\left( F(x)\right) ^{^{p}}F^{^{\prime }}(x)bx\leq \frac{1}{p+1}%
\left( \int_{a}^{b}s^{\frac{-1}{p}}(x)dx\right) ^{p}\int_{a}^{b}s(x)\left(
F^{^{\prime }}(x)\right) ^{p+1}dx.
\end{equation*}%
Substituting this inequality into (\ref{c1}) and proceeding as in the proof
of Theorem 2.11, we have the following results.

$\smallskip $

\textbf{Theorem 2.14.} \textit{\ Assume that }$r,$\textit{\ }$s$\textit{\ be
positive and continuous functions in }$(a,b)$\textit{\ and }$p$\textit{\ is
a positive integer. Then }%
\begin{equation*}
\int_{a}^{b}r(x)\left( \int_{a}^{x}f(t)dt\right) ^{p+1}dx\leq
\int_{a}^{b}s(x)\left( f(x)\right) ^{p+1}dx,
\end{equation*}

\textit{for all functions }$f\geq 0$, \textit{where }$C=\sup_{a<x<b}\left(
\int_{x}^{b}r(t)dt\right) \left( \int_{a}^{b}s^{\frac{-1}{p}}(x)dx\right)
^{p}.$

$\smallskip $

\textbf{Theorem 2.15.} \textit{Assume that }$r,$\textit{\ }$s$\textit{\ be
positive and continuous functions on }$(a,b)$\textit{\ and }$p$\textit{\ is
a positive integer. Then }%
\begin{equation*}
\int_{a}^{b}r(x)\left( \int_{x}^{b}f(t)dt\right) ^{p+1}\leq
C\int_{a}^{b}s(x)\left( f(x)\right) ^{p+1}dx,
\end{equation*}%
\textit{for all functions }$f\geq 0$, \textit{where }$C=\sup_{a<x<b}\left(
\int_{a}^{x}r(t)dt\right) \left( \int_{a}^{b}s^{\frac{-1}{p}}(x)dx\right)
^{p}.$

$\smallskip $

In the following we apply the Yang inequality \cite{Y1}%
\begin{equation}
\int_{a}^{b}\left\vert y(t)\right\vert ^{p}\left\vert y^{^{\prime
}}(t)\right\vert ^{q}dt\leq \frac{q}{p+q}(b-a)^{p}\int_{a}^{b}\left\vert
y^{^{\prime }}(t)\right\vert ^{p+q}dt.  \label{Y1}
\end{equation}%
where $y~$is an absolutely continuous function on $[a,b]$ with $x(a)=0$ (or $%
y(b)=0),\,\ p\geq 0$ and $q\geq 1$, to prove a new inequality of Hardy's
type. Applying the inequality (\ref{H0}) and the inequality (\ref{Y1}) on
the term 
\begin{equation*}
\int_{a}^{b}R(x,b)F^{p}(x)F^{^{\prime }}(x)dx,
\end{equation*}%
we have 
\begin{eqnarray*}
\int_{a}^{b}R(x,b)F^{p}(x)F^{^{\prime }}(x)dx &\leq &\left(
\int_{a}^{b}R^{p}(x,b)dx\right) ^{\frac{1}{p}}\left(
\int_{a}^{b}F^{pq}(x)\left( F^{^{\prime }}(x)\right) ^{q}dx\right) ^{\frac{1%
}{q}} \\
&\leq &\left( \frac{1}{p+1}\right) ^{\frac{1}{q}}(b-a)^{p}\left(
\int_{a}^{b}R^{p}(x,b)dx\right) ^{\frac{1}{p}} \\
&&\times \left( \int_{a}^{b}\left( F^{^{\prime }}(x)\right)
^{q(p+1)}dx\right) ^{\frac{1}{q}},
\end{eqnarray*}%
where $1/p+1/q=1.$ Using this inequality and proceeding as in the proof of
Theorem 2.11, we get the following results.

\smallskip

\textbf{Theorem 2.16.}\textit{\ Assume that }$r$\textit{\ be positive and
continuous function on }$(a,b)$\textit{\ and }$p>1$\textit{\ is a positive
integer. Then }%
\begin{equation*}
\int_{a}^{b}r(x)\left( \int_{a}^{x}f(t)dt\right) ^{p+1}dx\leq C\left(
\int_{a}^{b}\left( f(x)\right) ^{\frac{p(p+1)}{p-1}}dx\right) ^{\frac{p-1}{p}%
},
\end{equation*}%
\textit{for all functions }$f\geq 0,$ \textit{where }$C=\left( p+1\right) ^{%
\frac{1}{p}}(b-a)^{p}\left( \int_{a}^{b}R^{p}(x,b)dx\right) ^{\frac{1}{p}}.$

\smallskip

\textbf{Theorem 2.17.}\textit{\ Assume that }$r$\textit{\ be positive and
continuous function on }$(a,b)$\textit{\ and }$p>1$\textit{\ is a positive
integer. Then }%
\begin{equation*}
\int_{a}^{b}r(x)\left( \int_{x}^{b}f(t)dt\right) ^{p+1}dx\leq \left(
\int_{a}^{b}\left( f(x)\right) ^{\frac{p(p+1)}{p-1}}dx\right) ^{\frac{p-1}{p}%
},
\end{equation*}%
\textit{for all functions }$f\geq 0$, \textit{where }$C=\left( p+1\right) ^{%
\frac{1}{p}}(b-a)^{p}\left( \int_{a}^{b}R^{p}(a,x)dx\right) ^{\frac{1}{p}}.$

From Theorems 2.16, 2.17, we have the following inequality.

$\smallskip $

\textbf{Corollary 2.1.} \textit{Assume that }$r$\textit{\ be positive and
continuous function on }$(a,b)$\textit{\ and }$p>1$\textit{\ is a positive
integer. Then}%
\begin{equation*}
\left( \int_{a}^{b}r(x)\left( \int_{a}^{x}f(t)dt\right) ^{p+1}dx\right) ^{%
\frac{1}{p+1}}\leq C_{1}\left( \int_{a}^{b}\left( f(x)\right) ^{\frac{p(p+1)%
}{p-1}}dx\right) ^{\frac{p-1}{p(p+1)}},
\end{equation*}%
\textit{for all functions }$f\geq 0,$\textit{\ where }%
\begin{equation*}
C=\left( p+1\right) ^{\frac{1}{p(p+1)}}(b-a)^{\frac{p}{p+1}}\left(
\int_{a}^{b}R^{p}(x,b)dx\right) ^{\frac{1}{p(p+1)}},
\end{equation*}
\textit{and} \textit{\ }%
\begin{equation*}
\left( \int_{a}^{b}r(x)\left( \int_{x}^{b}f(t)dt\right) ^{p+1}dx\right) ^{%
\frac{1}{p+1}}\leq C_{2}\left( \int_{a}^{b}\left( f(x)\right) ^{\frac{p(p+1)%
}{p-1}}dx\right) ^{\frac{p-1}{p(p+1)}},
\end{equation*}%
\textit{for all functions }$f\geq 0,$\textit{\ where }%
\begin{equation*}
C=\left( p+1\right) ^{\frac{1}{p(p+1)}}(b-a)^{\frac{p}{p+1}}\left(
\int_{a}^{b}R^{p}(a,x)dx\right) ^{\frac{1}{p(p+1)}}.
\end{equation*}

By choosing $q+1=\frac{p(p+1)}{p-1}$, we have from Corollary 2.1, the
following inequalities which are of the form (\ref{h1}).

$\smallskip $

\textbf{Corollary 2.2.} \textit{Assume that }$r$\textit{\ be positive and
continuous function on }$(a,b)$\textit{\ and }$p>1$\textit{\ is a positive
integer. Then }%
\begin{equation*}
\left( \int_{a}^{b}r(x)\left( \int_{a}^{x}f(t)dt\right) ^{p+1}dx\right) ^{%
\frac{1}{p+1}}\leq C_{\ast }\left( \int_{a}^{b}\left( f(x)\right)
^{q+1}dx\right) ^{\frac{1}{q+1}},
\end{equation*}%
\textit{for all functions }$f\geq 0$, \textit{where }%
\begin{equation*}
C_{\ast }=\left( p+1\right) ^{\frac{1}{p(p+1)}}(b-a)^{\frac{p}{p+1}}\left(
\int_{a}^{b}R^{p}(x,b)dx\right) ^{\frac{1}{p(p+1)}},
\end{equation*}%
\textit{and} 
\begin{equation*}
\left( \int_{a}^{b}r(x)\left( \int_{x}^{b}f(t)dt\right) ^{p+1}dx\right) ^{%
\frac{1}{p+1}}\leq C^{\ast }\left( \int_{a}^{b}\left( f(x)\right)
^{q+1}dx\right) ^{\frac{1}{q+1}},
\end{equation*}%
\textit{for all functions }$f\geq 0$, \textit{where }%
\begin{equation*}
C^{\ast }=\left( p+1\right) ^{\frac{1}{p(p+1)}}(b-a)^{\frac{p}{p+1}}\left(
\int_{a}^{b}R^{p}(a,x)dx\right) ^{\frac{1}{p(p+1)}}.
\end{equation*}

\textit{\smallskip }

Yang \cite{Y2} extended the inequality (\ref{Y1}) and proved that if $r(t)$
is a positive bounded function and $y~$is an absolutely continuous on $[a,b]$
with $y(a)=0$ (or $y(b)=0)$, $p\geq 0,$ $q\geq 1$, then%
\begin{equation}
\int_{a}^{b}r(t)\left\vert y(t)\right\vert ^{p}\left\vert y^{^{\prime
}}(t)\right\vert ^{q}dt\leq \frac{q}{p+q}(b-a)^{p}\int_{a}^{b}r(t)\left\vert
y^{^{\prime }}(t)\right\vert ^{p+q}dt.  \label{Y2}
\end{equation}%
Applying the inequality (\ref{Y2}) on the term $%
\int_{a}^{b}R(x,b)F^{p}(x)F^{^{\prime }}(x)dx,$ we have 
\begin{eqnarray*}
\int_{a}^{b}R(x,b)F^{p}(x)F^{^{\prime }}(x)dx &\leq &\frac{1}{p+1}%
(b-a)^{p}\left( \int_{a}^{b}R(x,b)\left( F^{^{\prime }}(x)\right)
^{p+1}dx\right)  \\
&=&\frac{1}{p+1}(b-a)^{p}\left( \int_{a}^{b}R(x,b)\left( f(x)\right)
^{p+1}dx\right) .
\end{eqnarray*}%
Using this inequality and proceeding as in the proof of Theorem 2.11, we
obtain the following results.

$\smallskip $

\textbf{Theorem 2.18.}\textit{\ Assume that }$r$\textit{\ be positive and
continuous function on }$(a,b)$\textit{\ and }$p\geq 0$\textit{\ is a
positive integer. Then }%
\begin{equation*}
\int_{a}^{b}r(x)\left( \int_{a}^{x}f(t)dt\right) ^{p+1}dx\leq
(b-a)^{p}\int_{a}^{b}R(x,b)\left( f(x)\right) ^{p+1}dx,
\end{equation*}%
\textit{for all functions }$f\geq 0.$

\smallskip

\textbf{Theorem 2.19.}\textit{\ Assume that }$r$\textit{\ be positive and
continuous function on }$(a,b)$\textit{\ and }$p\geq 0$\textit{\ is a
positive integer. Then }%
\begin{equation*}
\int_{a}^{b}r(x)\left( \int_{x}^{b}f(t)dt\right) ^{p+1}dx\leq
(b-a)^{p}\int_{a}^{b}R(a,x)\left( f(x)\right) ^{p+1}dx,
\end{equation*}%
\textit{for all functions }$f\geq 0.$

$\smallskip $

In the following, we apply an inequality due to Boyd \cite{Boyd}\textbf{\ }%
and the H\"{o}lder inequality to obtain new inequalities. The Boyd
inequality states that: If $y\in C^{1}[a,$\ $b]$\ with $y(a)=0$\ (or $y(b)=0)
$, then%
\begin{equation}
\int_{a}^{b}\left\vert y(t)\right\vert ^{\nu }\left\vert y^{^{\prime
}}(t)\right\vert ^{\eta }dt\leq N(\nu ,\eta ,s)(b-a)^{\nu }\left(
\int_{a}^{b}\left\vert y^{^{\prime }}(t)\right\vert ^{s}dt\right) ^{\frac{%
\nu +\eta }{s}},  \label{Boyd}
\end{equation}%
where $\nu >0$, $s>1$, $0\leq \eta <s$, 
\begin{equation}
N(\nu ,\eta ,s):=\frac{\left( s-\eta \right) \nu ^{\nu }\sigma ^{\nu +\eta
-s}}{(s-1)(\nu +\eta )\left( I(\nu ,\eta ,s)\right) ^{\nu }},\text{ }\sigma
:=\left\{ \frac{\nu (s-1)+(s-\eta )}{(s-1)(\nu +\eta )}\right\} ^{\frac{1}{s}%
},  \label{n}
\end{equation}%
and%
\begin{equation*}
I(\nu ,\eta ,s):=\int_{0}^{1}\left\{ 1+\frac{s(\eta -1)}{s-\eta }t\right\}
^{-(\nu +\eta +s\nu )/s\nu }[1+(\eta -1)t]t^{1/\nu -1}dt.
\end{equation*}%
Applying the inequality (\ref{H0}) and the inequality (\ref{Boyd}) on the
term 
\begin{equation*}
\int_{a}^{b}R(x,b)F^{p}(x)F^{^{\prime }}(x)dx,
\end{equation*}%
we have 
\begin{eqnarray}
\int_{a}^{b}R(x,b)F^{p}(x)F^{^{\prime }}(x)dx &\leq &\left(
\int_{a}^{b}R^{p}(x,b)dx\right) ^{\frac{1}{p}}\left(
\int_{a}^{b}F^{pq}(x)\left( F^{^{\prime }}(x)\right) ^{q}dx\right) ^{\frac{1%
}{q}}  \notag \\
&\leq &N^{\frac{1}{q}}(pq,q,s)(b-a)^{p}\left(
\int_{a}^{b}R^{p}(x,b)dx\right) ^{\frac{1}{p}}  \notag \\
&&\times \left( \int_{a}^{b}\left( F^{^{\prime }}(x)\right) ^{s}dx\right) ^{%
\frac{p+1}{s}},  \label{2.20}
\end{eqnarray}%
where $1/p+1/q=1,$ and $N(pq,q,s)$ is determined from (\ref{n}) by putting $%
\nu =pq$ and $\eta =q$. Using the inequality (\ref{2.20}) and proceeding as
in the proof of Theorem 2.11, we have the following results.

$\smallskip $

\textbf{Theorem 2.20}. \textit{Assume that }$r$\textit{\ be a positive
functions in }$(a,b),$\textit{\ }$p>0,$\textit{\ }$s>1$, $0\leq q<s$ \textit{%
and }$1/p+1/q=1$\textit{. Then }%
\begin{equation*}
\int_{a}^{b}r(x)\left( \int_{a}^{x}f(t)dt\right) ^{p+1}dx\leq (p+1)C\left(
\int_{a}^{b}\left( f(x)\right) ^{s}dx\right) ^{\frac{(p+1)}{s}},
\end{equation*}%
\textit{for all functions }$f\geq 0$ \textit{where} $C=N^{\frac{1}{q}%
}(pq,q,s)(b-a)^{p}\left( \int_{a}^{b}R^{p}(x,b)dx\right) ^{\frac{1}{p}}.$

$\smallskip $

\textbf{Theorem 2.21}. \textit{Assume that }$r$\textit{\ be a positive and
continuous function on }$(a,b),$\textit{\ }$p>1,$\textit{\ }$1<q<s$ \textit{%
and }$1/p+1/q=1$\textit{. Then }%
\begin{equation*}
\int_{a}^{b}r(x)\left( \int_{x}^{b}f(t)dt\right) ^{p+1}dx\leq (p+1)C\left(
\int_{a}^{b}\left( f(x)\right) ^{s}dx\right) ^{\frac{p+1}{s}},
\end{equation*}%
\textit{for all functions }$f\geq 0$, \textit{where }$C=N^{\frac{1}{q}%
}(pq,q,s)(b-a)^{p}\left( \int_{a}^{b}R^{p}(a,x)dx\right) ^{\frac{1}{p}}.$

$\smallskip $

The inequality (\ref{Boyd}) has immediate application when $\eta =s$. In
this case the inequality (\ref{Boyd}) becomes%
\begin{equation}
\int_{a}^{b}\left\vert y(t)\right\vert ^{\nu }\left\vert y^{^{\prime
}}(t)\right\vert ^{\eta }dt\leq L(\nu ,\eta )(b-a)^{\nu }\left(
\int_{a}^{b}\left\vert y^{^{\prime }}(t)\right\vert ^{\eta }dt\right) ^{%
\frac{\nu +\eta }{\eta }},  \label{l0}
\end{equation}%
where 
\begin{equation}
L(\nu ,\eta ):=\frac{\eta \nu ^{\eta }}{\nu +\eta }\left( \frac{\nu }{\nu
+\eta }\right) ^{\frac{\nu }{\eta }}\left( \dfrac{\Gamma \left( \frac{\eta +1%
}{\eta }+\frac{1}{\nu }\right) }{\Gamma \left( \frac{\eta +1}{\eta }\right)
\Gamma \left( \frac{1}{\nu }\right) }\right) ^{\nu },  \label{l2}
\end{equation}%
and $\Gamma $ is the Gamma function. Applying the inequality (\ref{l0}) on
the term%
\begin{equation*}
\int_{a}^{b}F^{pq}(x)\left( F^{^{\prime }}(x)\right) ^{q}dx,
\end{equation*}%
we have 
\begin{equation}
\int_{a}^{b}F^{pq}(x)\left( F^{^{\prime }}(x)\right) ^{q}dx\leq
L(pq,q)(b-a)^{pq}\left( \int_{a}^{b}\left( F^{^{\prime }}(x)\right)
^{q}dx\right) ^{\frac{pq+q}{q}}  \label{z0}
\end{equation}%
where 
\begin{equation}
L(pq,q)=\frac{(pq)^{q}}{p+1}\left( \frac{p}{p+1}\right) ^{p}\left( \dfrac{%
\Gamma \left( \frac{q+1}{q}+\frac{1}{pq}\right) }{\Gamma \left( \frac{q+1}{q}%
\right) \Gamma \left( \frac{1}{pq}\right) }\right) .  \label{Lq}
\end{equation}%
Using (\ref{z0}), we see that 
\begin{eqnarray*}
\int_{a}^{b}R(x,b)F^{p}(x)F^{^{\prime }}(x)dx &\leq &\left(
\int_{a}^{b}R^{p}(x,b)dx\right) ^{\frac{1}{p}}\left(
\int_{a}^{b}F^{pq}(x)\left( F^{^{\prime }}(x)\right) ^{q}dx\right) ^{\frac{1%
}{q}} \\
&\leq &L^{\frac{1}{q}}(pq,q)(b-a)^{p}\left( \int_{a}^{b}R^{p}(x,b)dx\right)
^{\frac{1}{p}} \\
&&\times \left( \int_{a}^{b}\left( F^{^{\prime }}(x)\right) ^{q}dx\right)
^{p+1},
\end{eqnarray*}%
where $1/p+1/q=1.$ This gives us the following results.

$\smallskip $

\textbf{Theorem 2.22}. \textit{Assume that }$r$\textit{\ be positive and
continuous function in }$(a,b),$\textit{\ }$p,$ $q>1$ \textit{and }$1/p+1/q=1
$\textit{. Then }%
\begin{equation*}
\int_{a}^{b}r(x)\left( \int_{a}^{x}f(t)dt\right) ^{p+1}dx\leq (p+1)C\left(
\int_{a}^{b}\left( f(x)\right) ^{q}dx\right) ^{(p+1)},
\end{equation*}%
\textit{for all functions }$f\geq 0$\textit{, where }$C=L^{\frac{1}{q}%
}(pq,q)(b-a)^{p}\left( \int_{a}^{b}R^{p}(x,b)dx\right) ^{\frac{1}{p}}$ 
\textit{and} $L^{\frac{1}{q}}(pq,q)$\textit{\ is defined as in (\ref{Lq}).}

$\smallskip $

\textbf{Theorem 2.23}. \textit{Assume that }$r$\textit{\ be positive and
continuous function on }$(a,b),$\textit{\ }$p,$ $q>1$ \textit{and }$1/p+1/q=1
$\textit{. Then }%
\begin{equation*}
\int_{a}^{b}r(x)\left( \int_{x}^{b}f(t)dt\right) ^{p+1}dx\leq (p+1)C\left(
\int_{a}^{b}\left( f(x)\right) ^{q}dx\right) ^{(p+1)},
\end{equation*}%
\textit{for all functions }$f\geq 0$, \textit{where }$C=L^{\frac{1}{q}%
}(pq,q)(b-a)^{p}\left( \int_{a}^{b}R^{p}(a,x)dx\right) ^{\frac{1}{p}}$ 
\textit{and} $L^{\frac{1}{q}}(pq,q)$\textit{\ is defined as in (\ref{Lq}).}

\bigskip

Next in the following, we will apply some Opial type inequalities obtained
by Beesack and Das \cite{BD} (see also \cite{S9}) to prove some new
inequalities of Hardy's types with two different weighted functions. These
inequalities are presented in the following theorems.

\smallskip

\textbf{Theorem 2.24. }\textit{Let }$p,$\textit{\ }$q$\textit{\ be positive
real numbers such that} $p+q>1,$\textit{\ and let }$r,$ $s$\textit{\ be
positive continuous functions in }$(a,b)$\textit{\ such that }$%
\int_{a}^{b}s^{\frac{-1}{p+q-1}}(t)dt<\infty .$\textit{\ If }$%
y:[a,b]\rightarrow \mathbb{R}$\textit{\ is delta differentiable with }$%
y(a)=0,$ \textit{then }%
\begin{equation}
\int_{a}^{X}r(x)\left\vert y(x)\right\vert ^{p}\left\vert y^{^{\prime
}}(x)\right\vert ^{q}dx\leq K_{1}(a,X,p,q)\int_{a}^{X}s(x)\left\vert
y^{^{\prime }}(x)\right\vert ^{p+q}dx,  \label{z1}
\end{equation}%
\textit{where}%
\begin{eqnarray}
K_{1}(a,X,p,q) &=&\left( \frac{q}{p+q}\right) ^{\frac{q}{p+q}}  \notag \\
&&\times \left( \int_{a}^{X}(r(x))^{\frac{p+q}{p}}(s(x))^{-\frac{q}{p}%
}\left( \int_{a}^{x}s^{\frac{-1}{p+q-1}}(t)dt\right) ^{(p+q-1)}dx\right) ^{%
\frac{p}{p+q}}.  \label{z2}
\end{eqnarray}

\smallskip

\textbf{Theorem 2.25. }\textit{Assume that }$p,$\textit{\ }$q$\textit{\ be
positive real numbers such that }$p+q>1,$\textit{\ and let }$r,s$\textit{\
be nonnegative continuous functions in }$(a,b)$\textit{\ such that }$%
\int_{b}^{b}s^{\frac{-1}{p+q-1}}(t)dt<\infty .$\textit{\ If }$%
y:[a,b]\rightarrow \mathbb{R}$\textit{\ is delta differentiable with }$%
y(b)=0 $\textit{, then we have}%
\begin{equation}
\int_{a}^{b}r(x)\left\vert y(x)\right\vert ^{p}\left\vert y^{^{\prime
}}(x)\right\vert ^{q}dx\leq K_{2}(X,b,p,q)\int_{X}^{b}s(x)\left\vert
y^{^{\prime }}(x)\right\vert ^{p+q}dx,  \label{z4}
\end{equation}%
\textit{where}%
\begin{eqnarray}
K_{2}(X,b,p,q) &=&\left( \frac{q}{p+q}\right) ^{\frac{q}{p+q}}  \notag \\
&&\times \left( \int_{a}^{b}(r(x))^{\frac{p+q}{p}}(s(x))^{-\frac{q}{p}%
}\left( \int_{x}^{b}s^{\frac{-1}{p+q-1}}(t)dt\right) ^{(p+q-1)}dx\right) ^{%
\frac{p}{p+q}}.  \label{z9}
\end{eqnarray}%
In the following, we assume that there exists $h\in (a,b)$ which is the
unique solution of the equation%
\begin{equation*}
K(p,q)=K_{1}(a,h,p,q)=K_{2}(h,b,p,q)<\infty ,
\end{equation*}%
where $K_{1}(a,h,p,q)$ and $K_{2}(h,b,p,q)$ are defined as in Theorems 2.24
and 2.25. The combination of Theorems 2.24 and 2.5 gives the following
result.

\smallskip

\textbf{Theorem 2.26. }\textit{Let }$p,$\textit{\ }$q$\textit{\ be positive
real numbers such that }$pq>0$ and $p+q>1,$\textit{\ and let }$r,$ $s$%
\textit{\ be nonnegative continuous functions on }$(a,b)$\textit{\ such that 
}$\int_{a}^{b}s^{\frac{-1}{p+q-1}}(t)dt<\infty .$\textit{\ If }$%
y:[a,b]\rightarrow \mathbb{R}$\textit{\ is delta differentiable with }$%
y(a)=0=y(b)$\textit{, then we have}%
\begin{equation}
\int_{a}^{b}r(x)\left\vert y(x)\right\vert ^{p}\left\vert y^{\prime
}(x)\right\vert ^{q}dx\leq K(p,q)\int_{a}^{b}s(x)\left\vert y^{^{\prime
}}(x)\right\vert ^{p+q}dx.  \label{z13}
\end{equation}

For $r=s$ in (\ref{z1}), we obtain the following special case from Theorem
2.24.

\smallskip

\textbf{Corollary 2.3. }\textit{Let }$p,$\textit{\ }$q$\textit{\ be positive
real numbers such that }$p+q>1,$\textit{\ and let }$r$ \textit{be a
nonnegative continuous function in }$(a,b)$\textit{\ such that }$%
\int_{a}^{b}r^{\frac{-1}{p+q-1}}(t)dt<\infty .$\textit{\ If }$%
y:[a,b]\rightarrow \mathbb{R}$\textit{\ is delta differentiable with }$%
y(a)=0 $\textit{, then we have}%
\begin{equation}
\int_{a}^{X}r(x)\left\vert y(x)\right\vert ^{p}\left\vert y^{^{\prime
}}(x)\right\vert ^{q}dx\leq K_{1}^{\ast }(a,b,p,q)\int_{a}^{b}r(x)\left\vert
y^{^{\prime }}(x)\right\vert ^{p+q}dx,  \label{z34}
\end{equation}%
\textit{where}%
\begin{equation}
K_{1}^{\ast }(a,b,p,q)=\left( \frac{q}{p+q}\right) ^{\frac{q}{p+q}}\times
\left( \int_{a}^{b}r(x)\left( \int_{a}^{x}r^{\frac{-1}{p+q-1}}(t)dt\right)
^{(p+q-1)}dx\right) ^{\frac{p}{p+q}}.  \label{z35}
\end{equation}%
Now, we apply Theorem 2.24 to prove a new inequality of Hardy's type. Using
the fact that 
\begin{equation*}
\int_{a}^{b}R(x,b)F^{p}(x)F^{^{\prime }}(x)dx=\int_{a}^{b}\frac{R(x,b)}{r^{%
\frac{1}{q}}(x)}F^{p}(x)(r^{\frac{1}{q}}(x)F^{^{\prime }}(x))dx
\end{equation*}%
and applying the inequality (\ref{H0}), we have 
\begin{eqnarray}
\int_{a}^{b}\frac{R(x,b)}{r^{\frac{1}{q}}(x)}F^{p}(x)(r^{\frac{1}{q}%
}(x)F^{^{\prime }}(x))dx &\leq &\left( \int_{a}^{b}\frac{R^{p}(x,b)}{r^{%
\frac{1}{q}}(x)}dx\right) ^{\frac{1}{p}}  \notag \\
&&\times \left( \int_{a}^{b}r(x)F^{pq}(x)\left( F^{^{\prime }}(x)\right)
^{q}dx\right) ^{\frac{1}{q}}.  \label{a0}
\end{eqnarray}%
Applying the inequality (\ref{z1}) on the term 
\begin{equation*}
\int_{a}^{b}r(x)F^{pq}(x)\left( F^{^{\prime }}(x)\right) ^{q}dx,
\end{equation*}%
with $y=F$, \ and $p$ is replaced by $pq$, note that $F(a)=0$, we have that 
\begin{equation}
\int_{a}^{b}r(x)F^{pq}(x)\left( F^{^{\prime }}(x)\right) ^{q}dx\leq
K_{1}(a,b,pq,q)\int_{a}^{b}s(x)\left( F^{^{\prime }}(x)\right) ^{pq+q}dx,
\label{aq}
\end{equation}%
where 
\begin{eqnarray}
K_{1}(a,b,pq,q) &=&\left( \frac{q}{pq+q}\right) ^{\frac{q}{pq+q}}  \notag \\
&&\times \left( \int_{a}^{b}(r(x))^{\frac{pq+q}{p}}(s(x))^{-\frac{q}{pq}%
}\left( \int_{a}^{x}s^{\frac{-1}{pq+q-1}}(t)dt\right) ^{(pq+q-1)}dx\right) ^{%
\frac{pq}{pq+q}}.  \label{aq0}
\end{eqnarray}%
Substituting (\ref{aq}) into (\ref{a0}), we have that 
\begin{eqnarray*}
\int_{a}^{b}\frac{R(x,b)}{r^{\frac{1}{q}}(x)}F^{p}(x)(r^{\frac{1}{q}%
}(x)F^{^{\prime }}(x))dx &\leq &K_{1}^{\frac{1}{q}}(a,X,pq,q)\left(
\int_{a}^{b}\frac{R^{p}(x,b)}{r^{\frac{1}{q}}(x)}dx\right) ^{\frac{1}{p}} \\
&&\times \left( \int_{a}^{b}s(x)\left( F^{^{\prime }}(x)\right)
^{pq+q}dx\right) ^{\frac{1}{q}}.
\end{eqnarray*}%
This gives us the following result.

$\smallskip $

\textbf{Theorem 2.27.} \textit{Let }$p,$\textit{\ }$q$\textit{\ be positive
real numbers such that} $p,$ $q>1,$\textit{\ and let }$r,$ $s$\textit{\ be
nonnegative continuous functions on }$(a,b)$\textit{\ such that }$%
\int_{a}^{b}s^{\frac{-1}{p+q-1}}(t)dt<\infty .$ \textit{Then }%
\begin{eqnarray*}
\int_{a}^{b}r(x)\left( \int_{a}^{x}f(t)dt\right) ^{p+1}dx &\leq &(p+1)K_{1}^{%
\frac{1}{q}}(a,X,pq,q)\left( \int_{a}^{b}\frac{R^{p}(x,b)}{r^{\frac{1}{q}}(x)%
}dx\right) ^{\frac{1}{p}} \\
&&\left( \int_{a}^{b}s(x)\left( f(x)\right) ^{pq+q}dx\right) ^{\frac{1}{q}},
\end{eqnarray*}%
\textit{for all functions }$f\geq 0$, \textit{where }$K_{1}(a,b,pq,q)$%
\textit{\ is defined as in (\ref{aq0}).}

$\smallskip $

Applying Theorem 2.25 will give us the following result.

$\smallskip $

\textbf{Theorem 2.28.} \textit{Let }$p,$\textit{\ }$q$\textit{\ be positive
real numbers such that} $p,$ $q>1$ and $1/p+1/q=1,$\textit{\ and let }$r,$ $s
$\textit{\ be nonnegative continuous functions in }$(a,b)$\textit{\ such
that }$\int_{a}^{b}s^{\frac{-1}{p+q-1}}(t)dt<\infty .$ \textit{Then }%
\begin{eqnarray*}
\int_{a}^{b}r(x)\left( \int_{x}^{b}f(t)dt\right) ^{p+1}dx &\leq &(p+1)K_{2}^{%
\frac{1}{q}}(a,X,pq,q)\left( \int_{a}^{b}\frac{R^{p}(a,x)}{r^{\frac{1}{q}}(x)%
}dx\right) ^{\frac{1}{p}} \\
&&\left( \int_{a}^{b}s(x)\left( f(x)\right) ^{pq+q}dx\right) ^{\frac{1}{q}},
\end{eqnarray*}%
\textit{for all functions }$f\geq 0$, \textit{where }$K_{2}(a,b,pq,q)$%
\textit{\ is defined by }%
\begin{eqnarray*}
K_{2}(a,b,pq,q) &=&\left( \frac{q}{pq+q}\right) ^{\frac{q}{pq+q}} \\
&&\times \left( \int_{a}^{b}(r(x))^{\frac{pq+q}{pq}}(s(x))^{-\frac{q}{pq}%
}\left( \int_{x}^{b}s^{\frac{-1}{pq+q-1}}(t)dt\right) ^{(pq+q-1)}dx\right) ^{%
\frac{pq}{pq+q}}.
\end{eqnarray*}

\begin{remark}
One can apply Theorem 2.6 and Corollary 2.3 to obtain new inequalities of
Hardy's type. The details are left to the interested reader.
\end{remark}

In the following, we will apply a new Opial type inequality due to Beesack 
\cite{BSe3} to prove new inequalities of Hardy's type. The inequality due to
Beesack is given in the following theorem.

\bigskip

\textbf{Theorem 2.29.}\textit{\ Let }$r,$\textit{\ }$s$\textit{\ be
nonnegative, measurable functions on }$(\alpha ,\tau )$\textit{. Further
assume that }$k>1$\textit{, }$p>0$\textit{, }$0<q<k$\textit{, and let }$y(t)$%
\textit{\ be absolutely continuous in }$(\alpha ,\tau )$\textit{\ such that }%
$y(\alpha )=0.$\textit{\ Then }%
\begin{equation}
\int_{\alpha }^{\tau }r(t)\left\vert y(t)\right\vert ^{p}\left\vert
y^{^{\prime }}(t)\right\vert ^{q}dt\leq K_{1}(p,q,k)\left[ \int_{\alpha
}^{\tau }s(t)\left\vert y^{^{\prime }}(t)\right\vert ^{k}dt\right]
^{(p+q)/k},  \label{BS1}
\end{equation}%
\textit{where }%
\begin{eqnarray}
K_{1}(p,q,k) &=&\left( \frac{q}{q+p}\right) ^{\frac{q}{k}}  \notag \\
&&\times \left( \int_{\alpha }^{\tau }(r(y))^{\frac{k}{k-q}}(s(y))^{-\frac{q%
}{k-q}}\left( \int_{a}^{y}s^{\frac{-1}{k-1}}(t)dt\right)
^{p(k-1)/(k-q)}dy\right) ^{\frac{k-q}{k}}.  \label{ks1}
\end{eqnarray}%
\textit{If instead of }$(\alpha ,\tau )$\textit{\ is replaced by }$\mathit{(}%
\tau ,\beta )$\textit{\ and }$y(\alpha )=0$\textit{\ is replaced by }$%
y(\beta )=0$\textit{, then }%
\begin{equation}
\int_{\tau }^{\beta }r(t)\left\vert y(t)\right\vert ^{p}\left\vert
y^{^{\prime }}(t)\right\vert ^{q}dt\leq K_{2}(p,q,k)\left[ \int_{\tau
}^{\beta }s(t)\left\vert y^{^{\prime }}(t)\right\vert ^{k}dt\right]
^{(p+q)/k},  \label{Bs2}
\end{equation}%
\textit{where }%
\begin{eqnarray}
K_{2}(p,q,k) &=&\left( \frac{q}{q+p}\right) ^{\frac{q}{k}}  \notag \\
&&\times \left( \int_{\tau }^{\beta }(r(y))^{\frac{k}{k-q}}(s(y))^{-\frac{q}{%
k-q}}\left( \int_{y}^{\beta }s^{\frac{-1}{k-1}}(t)dt\right)
^{p(k-1)/(k-q)}dy\right) ^{\frac{k-q}{k}}.  \label{ks2}
\end{eqnarray}

Now, we apply inequality (\ref{BS1}) and (\ref{Bs2}) to obtain new
inequalities of Hardy's type. First we consider (\ref{BS1}). Using the fact
that 
\begin{equation*}
\int_{a}^{b}R(x,b)F^{p}(x)F^{^{\prime }}(x)dx=\int_{a}^{b}\frac{R(x,b)}{r^{%
\frac{1}{q}}(x)}F^{p}(x)(r^{\frac{1}{q}}(x)F^{^{\prime }}(x))dx,
\end{equation*}%
and applying the inequality (\ref{H0}), we have 
\begin{eqnarray}
\int_{a}^{b}\frac{R(x,b)}{r^{\frac{1}{q}}(x)}F^{p}(x)(r^{\frac{1}{q}%
}(x)F^{^{\prime }}(x))dx &\leq &\left( \int_{a}^{b}\frac{R^{p}(x,b)}{r^{%
\frac{1}{q}}(x)}dx\right) ^{\frac{1}{p}}  \notag \\
&&\times \left( \int_{a}^{b}r(x)F^{pq}(x)\left( F^{^{\prime }}(x)\right)
^{q}dx\right) ^{\frac{1}{q}}.  \label{cq}
\end{eqnarray}%
Applying (\ref{BS1}) on the term%
\begin{equation*}
\int_{a}^{b}r(x)F^{pq}(x)\left( F^{^{\prime }}(x)\right) ^{q}dx,
\end{equation*}%
we have (note that $F(a)=0),$ that%
\begin{equation}
\int_{a}^{b}r(x)F^{pq}(x)\left( F^{^{\prime }}(x)\right) ^{q}dx\leq
K_{1}(pq,q,k)\left[ \int_{a}^{b}s(x)(F^{^{\prime }}(x))^{k}dx\right]
^{(pq+q)/k},  \label{Bs3}
\end{equation}%
where 
\begin{eqnarray}
K_{1}(pq,q,k) &=&\left( \frac{1}{1+p}\right) ^{\frac{q}{k}}  \notag \\
&&\times \left( \int_{a}^{b}(r(x))^{\frac{k}{k-q}}(s(x))^{-\frac{q}{k-q}%
}\left( \int_{a}^{t}s^{\frac{-1}{k-1}}(t)dt\right) ^{pq(k-1)/(k-q)}dx\right)
^{\frac{k-q}{k}}.  \label{ks3}
\end{eqnarray}%
Substituting (\ref{Bs3}) into (\ref{cq}), we get that 
\begin{eqnarray*}
&&\int_{a}^{b}\frac{R(x,b)}{r^{\frac{1}{q}}(x)}F^{p}(x)(r^{\frac{1}{q}%
}(x)F^{^{\prime }}(x))dx \\
&\leq &K_{1}^{\frac{1}{q}}(pq,q,k)\left( \int_{a}^{b}\frac{R^{p}(x,b)}{r^{%
\frac{1}{q}}(x)}dx\right) ^{\frac{1}{p}}\times \left[ \int_{a}^{b}s(x)(F^{^{%
\prime }}(x))^{k}dx\right] ^{(p+1)/k}.
\end{eqnarray*}%
This gives us the following results.

\bigskip

\textbf{Theorem 2.30.} \textit{Let }$p,$\textit{\ }$q$\textit{\ be positive
real numbers such that} $p,$ $q>1,$\textit{\ and let }$r,$ $s$\textit{\ be
nonnegative continuous functions on }$(a,b).$ \textit{Then }%
\begin{eqnarray*}
\int_{a}^{b}r(x)\left( \int_{a}^{x}f(t)dt\right) ^{p+1}dx &\leq &(p+1)K_{1}^{%
\frac{1}{q}}(pq,q,q)\left( \int_{a}^{b}\frac{R^{p}(x,b)}{r^{\frac{1}{q}}(x)}%
dx\right) ^{\frac{1}{p}} \\
&&\times \left[ \int_{a}^{b}s(x)(f(x))^{k}dx\right] ^{(p+1)/k},
\end{eqnarray*}%
\textit{for all functions }$f\geq 0$, \textit{where }$K_{1}(pq,q,k)$\textit{%
\ is defined as in (\ref{ks3}).}

\bigskip

\textbf{Theorem 2.31.} \textit{Let }$p,$\textit{\ }$q$\textit{\ be positive
real numbers such that} $p,$ $q>1,$\textit{\ and let }$r,$ $s$\textit{\ be
nonnegative continuous functions on }$(a,b).$ \textit{Then }%
\begin{eqnarray*}
\int_{a}^{b}r(x)\left( \int_{x}^{b}f(t)dt\right) ^{p+1}dx &\leq &(p+1)K_{2}^{%
\frac{1}{q}}(pq,q,q)\left( \int_{a}^{b}\frac{R^{p}(a,x)}{r^{\frac{1}{q}}(x)}%
dx\right) ^{\frac{1}{p}} \\
&&\times \left[ \int_{a}^{b}s(x)(f(x))^{k}dx\right] ^{(p+1)/k},
\end{eqnarray*}%
\textit{for all functions }$f\geq 0$, \textit{where }$K_{2}(pq,q,k)$\textit{%
\ is defined by}%
\begin{eqnarray*}
K_{2}(pq,q,k) &=&\left( \frac{1}{1+p}\right) ^{\frac{q}{k}} \\
&&\times \left( \int_{a}^{b}(r(x))^{\frac{k}{k-q}}(s(x))^{-\frac{q}{k-q}%
}\left( \int_{t}^{b}s^{\frac{-1}{k-1}}(t)dt\right) ^{pq(k-1)/(k-q)}dx\right)
^{\frac{k-q}{k}}.
\end{eqnarray*}

\begin{remark}
The applications of the Opial inequality may be continue and one can obtain
new inequalities of Hardy's type. On the other hand one can obtain the
differential forms of the inequalities by replacing $\int_{a}^{x}f(t)dt$ by $%
g(x)$ on the right hand sides and replace $f(x)$ in the left hand sides by $%
g^{^{\prime }}(x).$ The details are left to the reader.
\end{remark}

\end{document}